\documentclass[10pt]{amsart}

\title[The limit of Betti numbers of a tower of finite covers]{The limit 
of $\mathbf{\IF_p}$-Betti numbers of a tower of 
finite covers with amenable fundamental groups}

\author{Peter Linnell}
\address{Department of Mathematics\\
Virginia Tech\\
Blacksburg, VA 24061-0123,
\\USA}
\email{plinnell@math.vt.edu}
\urladdr{http://www.math.vt.edu/people/plinnell/}

\author{Wolfgang L\"uck}
        \address{Westf\"alische Wilhelms-Universit\"at M\"unster\\
               Mathematisches Institut\\
               Einsteinstr.~62,
               D-48149 M\"unster, Germany}
         \email{lueck@math.uni-muenster.de}
         \urladdr{http://www.math.uni-muenster.de/u/lueck}

\author{Roman Sauer}
        \address{Westf\"alische Wilhelms-Universit\"at M\"unster\\
               Mathematisches Institut\\
               Einsteinstr.~62,
               D-48149 M\"unster, Germany}
         \email{sauerr@uni-muenster.de}
         \urladdr{http://wwwmath.uni-muenster.de/u/sauerr/}

\thanks{The authors thank the HIM at Bonn for its hospitality
		during the Trimester program ``Rigidity''  in the fall 2009 when this paper was written.
		This work was financially supported by the Leibniz-Preis of the second author.}
\date{March 1, 2010}

\usepackage{hyperref}
\usepackage[DIV10, headinclude, BCOR 10mm]{typearea}
\usepackage[initials,msc-links]{amsrefs}
\usepackage[initials]{amsrefs}
\usepackage{enumerate,amssymb}
\usepackage[arrow,matrix,tips]{xy}
  \SelectTips{eu}{10} \UseTips

\theoremstyle{plain}
\newtheorem{theorem}{Theorem}[section]
\newtheorem{lemma}[theorem]{Lemma}

\theoremstyle{definition}
\newtheorem{definition}[theorem]{Definition}
\newtheorem{example}[theorem]{Example}

\newtheorem{remark}[theorem]{Remark}

{\catcode`@=11\global\let\c@equation=\c@theorem}

\newcommand{\caln}{{\mathcal N}}

\newcommand{\calA}{{\mathcal A}}

\newcommand{\IC}{{\mathbb C}}

\newcommand{\IF}{{\mathbb F}}

\newcommand{\IN}{{\mathbb N}}

\newcommand{\IQ}{{\mathbb Q}}

\newcommand{\IZ}{{\mathbb Z}}

\newcommand{\aut}{\operatorname{aut}}

\newcommand{\charac}{\operatorname{char}}

\newcommand{\coker}{\operatorname{coker}}

\newcommand{\Elek}{\operatorname{Elek}}

\newcommand{\im}{\operatorname{im}}

\newcommand{\map}{\operatorname{map}}

\newcommand{\Ore}{\operatorname{Ore}}
\newcommand{\pr}{\operatorname{pr}}

\newcommand{\res}{\operatorname{res}}

\newcommand{\supp}{\operatorname{supp}}

\newcommand{\vdim}{\operatorname{vdim}}

\DeclareMathOperator{\diag}{diag}
\DeclareMathOperator{\Mat}{M}

\newcommand{\version}[1]                   
{\begin{center} last edited on #1\\
last compiled on \today\\
name of texfile: \jobname
\end{center}
}

\newcounter{commentcounter}

\begin{document}

\begin{abstract}
 We prove an analogue of the Approximation Theorem of $L^2$-Betti numbers by 
Betti numbers for arbitrary coefficient fields and virtually
torsionfree amenable groups. The limit of Betti numbers is identified
as the dimension of some module over the Ore localization of the group
ring. 
 \end{abstract}

\maketitle


\typeout{--------------------------------- Introduction  ------------------------------------}

\setcounter{section}{-1}
\section{Introduction}

A \emph{residual chain} of a group $G$ is a 
sequence $G=G_0\supset G_1\supset G_2\supset \cdots$ of normal subgroups 
of finite index such that $\bigcap_{i\ge 0}G_i=\{e\}$. 
The $n$-th $L^2$-Betti number of any 
finite free $G$-CW complex $X$ is the limit of the $n$-th Betti numbers of $G_i\backslash X$ 
normalized by the index $[G:G_i]$ for $i\to\infty$~\cite{Lueck(1994c)}. 
If we instead consider Betti numbers 
$b_n(G_i\backslash X;k)$ with respect to a field of characteristic $p>0$, the questions 
whether the limit exists, what it is, and whether it is independent of the residual chain 
are completely open for arbitrary residually finite $G$. 

For $G=\IZ^k$ 
and every field $k$ 
Elek showed that $\lim _{i\to\infty} b_n(G_i\backslash X;k)$ exists and expresses it in 
terms of the entropy of $G$-actions on the Pontrjagin duals of finitely generated 
$kG$-modules~\cite{Elek(2002)} -- his techniques play an important role in this 
paper~(see Section~\ref{sec:Eleks_dimension_function}). 
It was observed 
in~\cite{Abert-Jaikin-Zapirain-Nikolov(2007)}*{Theorem~17} 
that the mere convergence of the right hand side 
of~\ref{the:dim_approximation_over_fields:finitely_presented_modules} in 
Theorem~\ref{the:dim_approximation_over_fields} 
for every amenable $G$ and every field $k$ follows from a general 
convergence principle for subadditive functions on 
amenable groups~\cite{Lindenstrauss-Weiss(2000)}
and a theorem by Weiss~\cite{Weiss(2001)}. 

The main purpose of this paper is to determine 
the limit $\lim _{i\to\infty} b_n(G_i\backslash X;k)$ 
in algebraic terms 
for a large class of 
amenable groups including virtual torsionfree elementary amenable groups. 
This makes the limit computable 
by homological techniques; see e.g., the spectral sequence argument of 
Example~\ref{exa:s1_actions}. 

More precisely, the limit will be expressed in terms of the \emph{Ore dimension}. 
The group ring 
$kG$ of a torsionfree amenable group satisfying the zero-divisor conjecture 
fulfills the Ore condition with respect to the subset 
$S=kG-\{0\}$~\cite[Example~8.16 on page~324]{Lueck(2002)}; we will review 
the Ore localization in Subsection~\ref{sub:ore_localization}. 
The Ore localization $S^{-1}kG$ is a skew field containing $kG$.
Therefore the following definition makes sense: 

\begin{definition}[Ore dimension]
  \label{def:Ore_dimension}
  Let $G$ be a torsionfree amenable group such that $kG$ contains no 
  zero-divisors.  The \emph{Ore dimension} of a $kG$-module $M$ is defined by
\[\dim^{\Ore}_{kG}(M) = \dim_{S^{-1}kG}\bigl(S^{-1}kG \otimes_{kG} M\bigr).\]
\end{definition}

The following theorem is our main result; we will prove a more general version, including 
virtually torsionfree groups, in Section~\ref{sec:Extension_to_the_virtually_torsionfree_case}. 

\begin{theorem}\label{the:dim_approximation_over_fields}
Let $k$ be a field. Let $G$ be a torsionfree amenable group for which 
$kG$ has no zero-divisors\footnote{This assumption is
satisfied if $G$ is torsionfree elementary amenable. See Remark~\ref{rem:The_Ore_condition_for_group_rings}}. 
Let $(G_n)_{n\ge 0}$ be a residual chain of $G$. Then: 
\begin{enumerate}
%

\item  \label{the:dim_approximation_over_fields:finitely_presented_modules}
Consider a finitely presented $kG$-module $M$. Then
\begin{equation*}
\dim_{kG}^{\Ore}(M)
 = 
\lim_{n \to \infty} \frac{\dim_k\bigl(k \otimes_{kG_n}  M\bigr)}{[G:G_n]};
\end{equation*}

\item  \label{the:dim_approximation_over_fields:chain_complexes}
Consider a finite free $kG$-chain complex $C_*$. Then we get for all $i \ge 0$
\begin{equation*}
 \dim_{kG}^{\Ore}\bigl(H_i(C_*)\bigr)  
  = 
\lim_{n \to \infty} \frac{\dim_k\bigl(H_i(k \otimes_{kG_n} C_*)\bigr)}{[G:G_n]};
\end{equation*}

\item  \label{the:dim_approximation_over_fields:CW-complexes}
Let $X$ be a finite free $G$-$CW$-complex. Then we get for all $i \ge 0$
\begin{equation*}
\dim_{kG}^{\Ore} \bigl(H_i(X)\bigr)  
  = 
\lim_{n \to \infty} \frac{\dim_k\bigl(H_i(G_n \backslash X;k)\bigr)}{[G:G_n]}.
\end{equation*}
\end{enumerate}
\end{theorem}

%

\begin{remark}[Fields of characteristic zero] 
\label{rem:fields_of_characteristic_zero} 
Let $G$ be a group with a residual chain $(G_n)_{n \ge 0}$, and let $M$ be a 
finitely presented $kG$-module. Then 
the Approximation Theorem for $L^2$-Betti numbers says that 
\begin{equation}\label{eq:limit_for_modules}
\dim_{\caln(G)}\bigl(\caln(G) \otimes_{kG} M\bigr) = \lim_{n \to \infty} \frac{\dim_{k}\bigl(k \otimes_{kG_n} M\bigr)}{[G:G_n]}
\end{equation}
provided $k$ is an algebraic number field. Here 
$\caln(G)$ is the group von Neumann algebra, 
and $\dim_{\caln(G)}$ is the von Neumann dimension.
See~\cite{Lueck(1994c)} for $k = \mathbb{Q}$ and
\cite{Dodziuk-Linnell-Mathai-Schick_Yates(2003)} for the general case.

Let $k$ be a field of characteristic zero and let $u = \sum_{g \in G} x_g
\cdot g\in kG$ be an element. Let $F$ be the finitely generated
field extension of $\IQ$ given by $F = \IQ(x_g \mid g \in G) \subset k$.  
Then $u$ is already an element in $FG$. 
The field $F$ embeds into $\IC$: since $F$ is finitely generated, it is 
a finite algebraic extension of a transcendental extension $F'$ 
of~$\IQ$~\cite{Lang(2002_Algebra)}*{Theorem~1.1 on p.~356}, and 
$F'$ has finite transcendence degree over $\IQ$. Since the transcendence degree of 
$\IC$ over $\IQ$ is infinite, there exists an embedding $F'\hookrightarrow\IC$ induced 
by an injection of a transcendence basis of $F'/\IQ$ into a transcendence basis $\IC/\IQ$, 
which extends to $F\hookrightarrow\IC$ because $\IC$ is algebraically closed. This 
reduces the case of fields of characteristic zero to the case $k=\IC$. 
In~\cite{Elek(2006strong)} 
Elek proved~\eqref{eq:limit_for_modules} for amenable $G$ and 
$k=\IC$ (see also~\cite{Pape(2008)}). 

Moreover, 
if $G$ is a torsionfree amenable group such that $\IC G$ contains no zero-divisors
and $k$ is a field of characteristic zero, then
\[\dim_{\caln(G)}\bigl(\caln(G) \otimes_{kG} M\bigr) = \dim^{\Ore}_{kG}(M).\]
This follows from~\cite[Theorem~6.37 on page~259, Theorem~8.29 on page~330,
Lemma~10.16 on page~376, and Lemma~10.39 on page~388]{Lueck(2002)}. In particular, 
Theorem~\ref{the:dim_approximation_over_fields} follows for $k$ of characteristic zero.
So the interesting new case is the one of a field of prime characteristic. 
\end{remark}


\typeout{--------------------------------- Section 1  ----------------------------------------}

\section{Review of Ore localization and Elek's dimension function}
\label{sec:Review_of_Ore_localization_and_Eleks_dimension_function}


\subsection{Ore localization}\label{sub:ore_localization}
We review the Ore localization of rings. For proofs and more 
information the reader is referred to~\cite{Stenstroem(1975)}. 
Consider a torsionfree group $G$ and a field $k$. Let $S$ be the set of
non-zero-divisors of $k G$.  This is a multiplicatively closed subset of
$k G$ and contains the unit element of $k G$.  Suppose that $k G$ satisfies the
\emph{Kaplansky Conjecture} or \emph{zero-divisor conjecture}, 
i.e., $S=k G-\{0\}$.  Further assume that
$S$ satisfies the \emph{left Ore condition}, i.e., for $r \in k G$ and $s \in S$ there
exists $r' \in k G$ and $s' \in S$ with $s'r = r's$. Then we can consider the
\emph{Ore localization $S^{-1}k G$}.  Recall that every element in $S^{-1}kG$ is of the
form $s^{-1}\cdot r$ for $r \in k G$ and $s \in S$ and $s_0^{-1}\cdot r_0 =
s_1^{-1}\cdot r_1$ holds if and only there exists $u_0,u_1 \in R$ satisfying
$u_0 r_0 = u_1r_1$ and $u_0s_0 = u_1s_1 $.  Addition is given on representatives
by $s_0^{-1}r_0 + s_1^{-1}r_1 = t^{-1}(c_0r_0 + c_1r_1)$ for $t = c_0s_0 =
c_1s_1$.  Multiplication is given on representatives by $s_0^{-1}r_0 \cdot
s_1^{-1}r_1 = (ts_0)^{-1}cr_1$, where $cs_1 = t r_0$. The zero element is
$e^{-1}\cdot 0$ and the unit element is $e^{-1}\cdot e$.
The Ore localization $S^{-1}kG$
is a skew field and the canonical map $k G \to S^{-1}kG$ sending $r$ to
$e^{-1} \cdot r$
is injective. The functor $S^{-1}kG \otimes_{kG} -$ is exact. 

\begin{remark}[The Ore condition for group rings]
\label{rem:The_Ore_condition_for_group_rings}
If a torsionfree amenable group $G$ satisfies
the Kaplansky Conjecture, i.e., $kG$ contains no zero-divisor,
then for $S = kG-\{0\}$ the Ore localization $S^{-1}kG$ exists and is a 
skew field~\cite{Lueck(2002)}*{Example~8.16 on page~324}.
Every torsionfree elementary amenable group satisfies the assumptions above
for all fields $k$~\citelist{\cite{Kropholler-Linnell-Moody(1988)}*{Theorem~1.2}
\cite{Linnell(2006)}*{Theorem~2.3}}.  If the group $G$ contains the free group
of rank two as subgroup, then the Ore condition is never satisfied for 
$kG$~\cite{Linnell(2006)}*{Proposition~2.2}.
\end{remark}

From the previous remark and the discussion above we obtain: 

\begin{theorem}\label{the:Ore-dimension}
  Let $G$ be a torsionfree amenable group such that $kG$ contains no 
  zero-divisors.  Then the Ore dimension $\dim^{\Ore}_{kG}$ has the following
  properties:
  \begin{enumerate}
  \item $\dim^{\Ore}_{kG}(kG) = 1$;
  \item For any short exact sequence of $kG$-modules $0 \to M_0 \to M_1 \to
    M_2\to 0$ we get
$$\dim^{\Ore}_{kG}(M_1) = \dim^{\Ore}_{kG}(M_0) + \dim^{\Ore}_{kG}(M_2).$$
\end{enumerate}
\end{theorem}


\subsection{Crossed products, Goldie rings, and the generalized Ore localization }
Throughout, let $G$ be a group, let $k$ be a skew field. 

Let $R$ be a ring. The notion of crossed product 
generalizes the one of group ring. A \emph{crossed product} 
$R\ast G=R\ast_{c,\tau} G$ is determined by maps $c:G\to\aut(R)$ and 
$\tau:G\times G\to R^\times$ such that, roughly speaking, $c$ is a homomorphism up 
to the $2$-cocycle $\tau$. We refer to the 
survey~\cite{Lueck(2002)}*{10.3.2 on p.~398} for details. If $G$ is an extension 
of $H$ by $Q$, then the group ring $kG$ is isomorphic to a crossed product 
$kH\ast Q$. Some results in this paper are formulated for crossed products, although 
we only need the case of group rings for 
Theorem~\ref{the:dim_approximation_over_fields}. So the reader may think of 
group rings most of the time. However, crossed products show up naturally, e.g., 
in proving that the virtual Ore dimension~\eqref{eq:virtual_Ore_dimension} is 
well defined. 

We recall the following definition. 

\begin{definition}\label{def:goldie_ring}
	A ring $R$ is \emph{left Goldie} if there exists $d\in \mathbb{N}$ such
	that every direct sum of nonzero left ideals of $R$ has at most
	$d$ summands and the left annihilators $a(x)=\{r\in R;~rx=0\}$, $x\in R$, 
	satisfy the maximum condition for ascending chains. A ring $R$ is \emph{prime}
	if for any two ideals $A,B$ in $R$, $AB=0$ implies $A=0$ or $B=0$. 
\end{definition}

The subgroup of $G$ generated by its finite normal subgroups will be denoted by 
$\Delta^+(G)$.  Then $\Delta^+(G)$ is also the set of elements of
finite order which have only finitely many conjugates.
We need the following three results: 

\begin{lemma}[\cite{Passman(1986)}*{Corollary 5 of Lecture 4}]\label{lem:delta_implies_prime}
If $\Delta^+(G) = 1$, then $k*G$ is prime. 
\end{lemma}

\begin{theorem}[\cite{Passman(1977)}*{Theorem~4.10 on p.~456}]\label{thm:prime_Goldie_matrix}
	The set of non-zero-divisors in a prime left Goldie ring satisfies the Ore 
	condition. The Ore localization $S^{-1}R$ is isomorphic to 
	$\Mat_d(D)$ for some $d \in \mathbb{N}$ and
	some skew field $D$. 
\end{theorem}
 
\begin{theorem}\label{thm:elementary_amenable_goldie}
	If $G$ is amenable and $k*G$ is a
	domain, then $k*G$ is a prime left Goldie ring. 
	If $G$ is an elementary amenable group such that the orders
of the finite 
	subgroups are bounded, then $k*G$ is left Goldie. 
\end{theorem}

\begin{proof}
	If $G$ is amenable and $k\ast G$ is a domain, then $k\ast G$ satisfies the 
	Ore condition~\cite{Dodziuk-Linnell-Mathai-Schick_Yates(2003)}*{Theorem~6.3}, thus 
	its Ore localization with respect to $S=k\ast G-\{0\}$ 
	is a skew field. By~\cite{Passman(1977)}*{Theorem~4.10 on p.~456} 
	$k\ast G$ is a prime left Goldie ring. The second assertion is taken 
	from~\cite{Kropholler-Linnell-Moody(1988)}*{Proposition~4.2}. 
\end{proof}

Next we extend the definition of Ore dimension to prime left Goldie rings. 
Let $R$ be such a ring. 
The functor $S^{-1}R \otimes_R -$ will still be
exact~\cite{Stenstroem(1975)}*{Proposition~II.1.4 on page~51}. 
If $M$ is a left $R$-module, then $S^{-1}R \otimes_R M$ will
be a direct sum of $n$ irreducible $S^{-1}R$-modules for some
non-negative integer $n$, and then the \emph{(generalized) Ore dimension} 
of $M$ is defined as 
\[\dim_R^{\Ore}(M) = \frac{n}{d}.\]
Since $S^{-1}R\cong\Mat_d(D)$ (Theorem~\ref{thm:prime_Goldie_matrix}) and 
$\Mat_d(D)$ decomposes into $d$ copies of the irreducible module $D^d$, we have 
$\dim_R^{\Ore}(R)=1$.


\subsection{Elek's dimension function}
\label{sec:Eleks_dimension_function}

Throughout this subsection let $G$ be a finitely generated amenable group. 
We review Elek's definition~\cite{Elek(2003c)} 
of a dimension function $\dim^{\Elek}_{kG}$ for finitely generated 
$kG$-modules. 

Fix a finite set of generators and equip $G$ with the associated
word metric $d_G$. A \emph{F\o lner sequence} $(F_n)_{n \ge 0}$ is a 
sequence of finite subsets of $G$ such that for 
any fixed $R > 0$
we have
\[\lim_{n \to \infty} \frac{|\partial_RF_n|}{|F_n|} = 0,\]
where $\partial_RF_n=\{g \in G \mid d(g,F_k) \le R \; \text{and}\; d(g,G\setminus F_k) \le R\}$.

Let $k$ be an arbitrary skew field endowed with the discrete topology
and let $\mathbb{N}$ denote the positive integers $\{1,2,\dots\}$. Let $n\in\IN$. We equip 
the  space 
of functions $\map(G, k^n)=\prod_{g\in G}k^n$ with the product topology, which is the same as the 
topology of pointwise convergence. The natural right $G$-action on 
$\map(G, k^n)$ is defined by 
\[
	(\phi g)(x)=\phi(xg^{-1})\text{ for $g,x\in G,\phi\in\map(G,k^n)$.}
\]
Also $\map(G,k^n)$ is a right $k$-vector space by defining $(\phi
k)(x) = \phi(x)k$.
For any subset $S\subset G$ and any subset $W\subset\map(G,k^n)$ let 
\[
	W\vert_S=\{f:S\to k^n \mid \exists g\in W \;\text{with} \; g\vert_S=f\}. 
\]
A right $k$-linear subspace $V\subset\map(G,k^n)$ is called \emph{invariant} if $V$ is closed and 
invariant under the right $G$-action. 

Elek defines the \emph{average dimension} 
$\dim_G^\calA(V)$ of an invariant subspace $V$ by choosing a F\o lner sequence 
$(F_n)_{n\in\IN}$ of $G$ and setting 
\begin{equation}\label{eq:average_dimension}
	\dim_G^\calA(V)=\limsup_{n\to\infty}\frac{\dim_k \bigl(V\vert_{F_n}\bigr)}{|F_n|}. 
\end{equation}

\begin{theorem}[\cite{Elek(2003c)}*{Prop.~7.2 and Prop.~9.2}]
	The sequence in~\eqref{eq:average_dimension} converges and its limit 
	$\dim_G^\calA(V)$ is independent of the choice of the F\o lner 
	sequence. 
\end{theorem}

\begin{remark}
     Elek actually defines $\dim_G^\calA(V)$ using F\o lner exhaustions,
     i.e.~increasing F\o lner sequences $(F_{n \in \mathbb{N}})$ with
     $\bigcup_{n\in\IN} F_n =G$.  This makes no difference since the
     existence of the limit of $(\dim_k \bigl(V\vert_{F_n}\bigr)/|F_n|)_{n\in\IN}$ 
	 for arbitrary F\o lner sequences (and thus its independence of the choice) follows 
	 from~\cite{Lindenstrauss-Weiss(2000)}*{Theorem~6.1}. 
\end{remark}

Let $M$ be a finitely generated left $kG$-module. The $k$-dual 
$M^\ast=\hom_k(M,k)$ (where $M$ and $k$ are viewed as left
$k$-modules, and $(\phi a)m = \phi(am)$ for $\phi \in M^*$, $a \in k$
and $m\in M$)
carries the natural right $G$-action $(\phi g)(m)=\phi(g m)$. 
The dual of the free left $kG$-module $kG^n$ is canonically isomorphic
to $\map(G, k^n)$. 
Any left $kG$-surjection $f\colon kG^n\twoheadrightarrow M$ induces a
right $kG$-injection
$f^\ast\colon M^\ast\to\map(G, k^n)$ such that
$\im(f^\ast)$ is a $G$-invariant $k$-subspace. 

\begin{definition}[Elek's dimension function]
 \label{def:elek_dimension_function} 
Let $M$ be a  finitely generated left $kG$-module. Its \emph{dimension in the sense of Elek} is 
defined by choosing a left $kG$-surjection $f\colon kG^n\twoheadrightarrow M$ and 
setting 
	\begin{equation}\label{eq:elek_dimension_and_average_dimension}
	\dim^{\Elek}_{kG}(M) =	\dim_G^\calA\bigl(\im(f^\ast)\bigr). 
	\end{equation}
\end{definition}

\begin{theorem}[Main properties of Elek's dimension function]\label{the:Eleks-dimension}
Let $G$ be a finitely generated amenable group. The 
definition~\eqref{eq:elek_dimension_and_average_dimension} of $\dim^{\Elek}_{kG}(M)$ is 
independent of the choice of the surjection $f$, and $\dim^{\Elek}_{kG}$
has the following properties:
\begin{enumerate}
  \item $\dim^{\Elek}_{kG}(kG) = 1$;
  \item For  any short exact sequence of finitely generated
$kG$-modules $0 \to M_0 \to M_1 \to M_2 \to 0$ we get
$$\dim^{\Elek}_{kG}(M_1) = \dim^{\Elek}_{kG}(M_0) + \dim^{\Elek}_{kG}(M_2);$$
  \item If the finitely generated $kG$-module $M$ satisfies $\dim^{\Elek}_{kG}(M) = 0$, then every quotient
  module $Q$ of $M$ satisfies $\dim^{\Elek}_{kG}(Q) = 0$.
\end{enumerate}
\end{theorem}

\begin{proof} The first two assertions are proved in~\cite[Theorem~1]{Elek(2003c)}.
Notice that the third condition does not necessarily follow from additivity
since the kernel of the epimorphism $M \to Q$ may not be finitely generated.
But the third statement is a direct consequence of the definition of Elek's dimension.
\end{proof}

\begin{remark}[The dual of finitely generated $kG$-modules]\label{rem:Pontrjagin_duality}
Identify the left $kG$-module $kG^n$ with the finitely supported
functions in $\map(G,k^n)$.  Here we view $\map(G,k^n)$ as a left
$k$-vector space by $(af)(g) = af(g)$, and the left $G$-action is
given by $(hf)(g) = f(h^{-1}g)$ for $h,g \in G$ and $a \in k$.
    Let 
	\begin{equation*}
	\langle\_,\_\rangle\colon kG^n \times \map(G, k^n)
	\end{equation*}
	be the canonical pairing (evaluation) of $kG^n$ and its dual
$\map(G, k^n)$.
	If we view an element $f\in kG^n$ as a finitely supported function 
	$G\to k^n$ (in $\map(G,k^n)$), then the pairing of $f\in kG^n$
with $l\in \map(G, k^n)$
is given by 
	\begin{equation*}
		\langle f,l\rangle=\sum_{g\in G}(f(g),l(g)),
	\end{equation*}
	where $(\_,\_)$ denotes the standard inner product in $k^n$. 
	For a subset $W\subset kG^n$ let 
	\[
		W^\perp = \bigl\{f\in\map(G,k^n) \mid \langle x, f\rangle=0~\forall x\in W\bigr\}. 
	\]
	If $M$ is a finitely generated $kG$-module and $f\colon kG^n\twoheadrightarrow M$ is 
	a left $kG$-surjection, then $f^* \colon M^* \hookrightarrow
\map(G,k^n)$ is a right $kG$-injection and
	\[
		\im(f^\ast)=\ker(f)^\perp\subseteq\map(G,k^n). 
	\]
\end{remark}


\typeout{--------------------------------- Section 2  ----------------------------------------}

\section{Approximation for finitely presented $kG$-modules for Elek's
  dimension function}
\label{sec:Approximation_for_finitely_presented_kG-modules_for_Eleks_dimension_function}

The main result of this section is: 

\begin{theorem} \label{the:Approximation_for_finitely_presented_kG-modules_for_Eleks_dimension_function}
  Let $G$ be a finitely generated amenable group. Consider a sequence
  of normal subgroups of finite index
$$G= G_0 \supseteq G_1 \supseteq G_2 \supseteq \cdots$$
such that $\bigcap_{n \ge 0} G_n = \{1\}$. Then every finitely
presented $kG$-module $M$ satisfies
\[
\dim_{kG}^{\Elek}(M)= \lim_{n \to \infty} \frac{\dim_k\bigl(k
  \otimes_{kG_n} M\bigr)}{[G:G_n]}.\]
\end{theorem}

Its proof needs some preparation.

Throughout, let $G$ be a finitely generated amenable group.  For
any subset $S\subset G$ let $k[S]$ be the $k$-subspace of $kG$
generated by $S\subset kG$. Let $j[S]\colon k[S]\to k[G]$ be the
inclusion and $\pr[S]\colon kG\to k[S]$ be the projection given by
\[ \pr[S](g)=\begin{cases}
  g & \text{if $g\in S$;}\\
  0 & \text{if $g\in G\backslash S$.}
\end{cases}\]

\begin{theorem}\label{the:equality_elek_and_approximative_dimension_of_finitely_presented_modules}
  Let $G$ be a finitely generated amenable group.  Let $M$ be a
  finitely presented left $kG$-module $M$ with a presentation
  $kG^r\xrightarrow{f}kG^s\xrightarrow{p} M\to 0$.  For every subset
  $S\subset G$ we define
  \[
  M[S]=\coker\Bigl(\pr[S]\circ f\circ j[S]\colon k[S]^r\to
  k[S]^s\Bigr).
  \]
  Let $(F_n)_{n \ge 0}$ be a F\o lner sequence of $G$. Then
  \[
  \dim^{\Elek}_{kG}(M)=\lim_{n\to\infty}\frac{\dim_k\left(M[F_n]\right)}{|F_n|}.
  \]
\end{theorem}

\begin{proof}
  The map $f$ is given by right multiplication with a matrix $A\in
  M_{r,s}(kG)$.  Viewing $A$ as a map $G\to k^{r\times s}$ it is clear
  what we mean by the support $\supp(A)$ of $A$. Let $R >0$ be the
  diameter of $\supp(A)\cup\supp(A)^{-1}$.  Since
  \[\lim_{n \to \infty} \frac{\bigl|\partial_R
    F_n\bigr|}{\bigl|F_n\bigr|} = 0,\] it is enough to show that for
  every $n\ge 1$
  \begin{equation}\label{eq:estimate_1}
    \bigl|\dim_k(M[F_n])-\dim_k\bigl(\im(p^\ast)\vert_{F_n}\bigr)\bigr| \le s\cdot \bigl|\partial_R F_n\bigr|. 
  \end{equation}
	
  For the definition of inner products $(\_,\_)$ and
  $\langle\_,\_\rangle$ we refer to Remark~\ref{rem:Pontrjagin_duality}.  Define the following $k$-linear subspaces of
  $\map(F_n,k^s)$:
  \begin{align*}
    W_n &= \bigl\{\phi\colon F_n\to k^s\mid
    \langle\pr_n\circ f \circ j_n(x),\phi\rangle=0~\forall x\in k[F_n]^r\bigr\};\\
    V_n &= \bigl\{\phi\colon F_n\to k^s \mid \exists \bar\phi:G\to k^s
    \; \text{satisfying}\;
    \bar\phi\vert_{F_n}=\phi, \langle f(y), \bar\phi\rangle=0~\forall y\in kG^r\bigr\};\\
    Z_n &=\bigl\{\phi\colon F_n\to k^s \mid \phi\vert_{\partial_R
      F_n}=0\bigr\}.
  \end{align*}

  Since $\dim_k(M[F_n]) = \dim_k(W_n)$ and $\dim_k\bigl(\im(p^\ast)\vert_{F_n}\bigr) = \dim_k(V_n)$, 
  the desired estimate~\eqref{eq:estimate_1} is equivalent to 
  \begin{equation}\label{eq:estimate_2}
    \bigl|\dim_k(W_n)-\dim_k(V_n)\bigr| \le s\cdot \bigl|\partial_R F_n\bigr|. 
  \end{equation}

  By additivity of $\dim_k$ we obtain that
  \begin{align*}
    \dim_k (W_n\cap Z_n)&\ge
    \dim_k(W_n)-\dim_k(\map(F_n,k^s))+\dim_k(Z_n)\\
    &\ge\dim_k(W_n)-s \cdot |F_n| + s \cdot (|F_n|-|\partial_R F_n|)\\
    &=\dim_k(W_n)-s \cdot |\partial_R F_n|.
  \end{align*}
  Similarly, we get
  \[
  \dim_k (V_n\cap Z_n)\ge \dim_k(V_n)-s \cdot |\partial_R F_n|.
  \]
  To prove~\eqref{eq:estimate_2} it hence suffices to show that
  \begin{align}
    W_n\cap Z_n&\subset V_n; \label{eq:first_inclusion}\\
    V_n\cap Z_n&\subset W_n. \label{eq:second_inclusion}
  \end{align}
  Let $\phi\in W_n\cap Z_n$. Extend $\phi$ by zero to a function
  $\bar\phi:G\to k^s$. Let $y\in kG^r$. Then we can decompose $y$ as
  $y=y_0+y_1$ with $\supp(y_0)\subset F_n$ and $\supp(y_1)\subset
  G\backslash F_n$. By definition of the radius $R$ it is clear that
  $\supp(f(y_1))\subset G\backslash F_n\cup\partial_RF_n$.  Because of
  $\phi\in Z_n$ we have $\langle f(y_1),\bar\phi\rangle=0$. The fact
  that $\phi\in W_n$ implies that
  \[\langle f(y_0), \bar\phi\rangle=\langle \pr_n\circ f\circ
  j_n(y_0), \phi\rangle=0.\] So we obtain that $\langle
  f(y),\bar\phi\rangle=0$, meaning that $\phi\in V_n$. The proof
  of~\eqref{eq:second_inclusion} is similar.
\end{proof}

The following theorem is due to Weiss. Its proof can be
found in~\cite[Proposition~5.5]{Deninger-Schmidt(2007)}.

\begin{theorem}[Weiss]\label{the:Weiss_Theorem}
  Let $G$ be a countable amenable group. Let $G_n\subset G$, $n\ge 1$,
  be a sequence of normal subgroups of finite index with
  $\bigcap_{n\ge 1}G_n = \{1\}$.  Then there exists, for every $R\ge
  1$ and every $\epsilon > 0$, an integer $M = M(R,\epsilon)\ge 1$
  such that for $n \ge M$ there is a fundamental domain $Q_n\subset G$
  of the coset space $G/G_n$ such that
  \[
  \frac{\bigl|\partial_R Q_n\bigr|}{\bigl|Q_n\bigr|} < \epsilon.
  \]
\end{theorem}

Now we are ready to prove
Theorem~\ref{the:Approximation_for_finitely_presented_kG-modules_for_Eleks_dimension_function}
\begin{proof}[Proof of
  Theorem~\ref{the:Approximation_for_finitely_presented_kG-modules_for_Eleks_dimension_function}]
  According to Theorem~\ref{the:Weiss_Theorem} let $(Q_n)_{n \ge 0}$
  be a F\o lner sequence of $G$ such that $Q_n\subset G$ is a
  fundamental domain for $G/G_n$.  Choose a finite presentation of
  $M$:
  \[
  kG^r\xrightarrow{f}kG^s\to M\to 0.
  \]
  Let $f_n=k[G/G_n]\otimes_{kG}f$. By right-exactness of tensor
  products we have the exact sequence
  \[
  k[G/G_n]^r\xrightarrow{f_n}k[G/G_n]^s\to k[G/G_n]\otimes_{kG} M\to
  0.
  \]
  The natural map $Q_n\subset G\rightarrow G/G_n$ induces an
  isomorphism $j_n\colon k[Q_n]\to k[G/G_n]$ of $k$-vector spaces.  The map
  $f$ is given by right multiplication $f=R_A$ with a matrix $A\in
  M_{r,s}(kG)$.  Viewing $A$ as a map $G\to k^{r\times s}$ let
  $\supp(A)$ be the support of $A$.  Let $R > 0$ be the diameter of
  $\supp(A)\cup\supp(A)^{-1}$ (with respect to the fixed word metric
  on $G$).  Then $f$ restricts to a map
  \[
  f\vert_{Q_n\setminus \partial_RQ_n}\colon
  k[Q_n\setminus \partial_RQ_n]^r\to k[Q_n]^s.
  \]
  Hence there is precisely one $k$-linear map $g$ for which the
  following diagram of $k$-vector spaces commutes:
  \begin{equation}\label{eq:approximation_diagram}
  \xymatrix@C=16mm{k[G/G_n]^r\ar[r]^{f_n}& k[G/G_n]^s\ar[r]&k[G/G_n]\otimes_{kG}M\ar[r]&0\\
    k[Q_n\setminus \partial_RQ_n]^r\ar[u]^{j_n\vert_{Q_n\setminus \partial_RQ_n}}
   \ar[r]^-{f\vert_{Q_n\setminus \partial_RQ_n}} &
    k[Q_n]^s\ar[u]^{j_n}_{\cong} \ar[r]^{\pr}&\coker(
    f\vert_{Q_n\setminus \partial_RQ_n})\ar@{->>}[u]^g\ar[r]&0}
  \end{equation}
  One easily verifies that $g$ is surjective and that
  \[
  \ker(g)\subset\im\bigl(\pr\circ j_n^{-1}\circ f_n\colon
  k[G/G_n]^r\to \coker( f\vert_{Q_n\setminus \partial_RQ_n})\bigr).
  \]
  The map $\pr\circ j_n^{-1}\circ f_n$ descends to a map
  \[
  \pr\circ j_n^{-1}\circ f_n\colon
  \coker(j_n\vert_{Q_n\setminus \partial_RQ_n})\to \coker(
  f\vert_{Q_n\setminus \partial_RQ_n}).
  \]
Note that
  \[
  \dim_k\bigl(\coker(j_n\vert_{Q_n\setminus \partial_RQ_n})\bigr)=
  r \cdot |\partial_R Q_n|.
  \]
  Thus,
  \begin{equation*}
    \dim_k(\coker(	f\vert_{Q_n\setminus \partial_RQ_n}))-\dim_k\left(k[G/G_n]\otimes_{kG}M\right)
    =\dim_k \ker(g)\le r\cdot \bigl|\partial_R Q_n\bigr|.
  \end{equation*}
  By replacing the upper row in diagram~\eqref{eq:approximation_diagram} by 
\[k[Q_n]^r\xrightarrow{\pr[Q_n]\circ f\circ j[Q_n]}k[Q_n]^s\to M[Q_n]\to 0\]
 and essentially running the same argument as before 
  we obtain that
  \begin{equation*}
    \dim_k(\coker(	f\vert_{Q_n\setminus
\partial_RQ_n}))-\dim_k(\coker(M[Q_n]))\le r\cdot  \bigl|\partial_R Q_n\bigr|. 
  \end{equation*}
  Since
  \[
  \frac{|\partial_R Q_n|}{[G/G_n]}=\frac{ |\partial_R Q_n|}{
    |Q_n|}\xrightarrow{n\to\infty} 0
  \]
  we get that
  \[
  \lim_{n\to\infty}\frac{\dim_k( k[G/G_n]\otimes_{kG}M)}{[G:G_n]}
  \]
  exists if and only if
  \begin{equation*}
    \lim_{n\to\infty}\frac{\dim_k( M[Q_n])}{|Q_n|}
  \end{equation*}
  exists, and in this case they are equal. Now the assertion follows
  from
  Theorem~\ref{the:equality_elek_and_approximative_dimension_of_finitely_presented_modules}.
\end{proof}


\typeout{--------------------------------- Section 3  ----------------------------------------}

\section{Comparing dimensions}
\label{sec:Comparing_dimensions}

The main result of this section is: 

\begin{theorem}[Comparing dimensions]
\label{the:Comparing_dimensions}
Let $G$ be a group, let
$k$ be a skew field, and let $k*G$ be a crossed product which is
prime left Goldie.
Let $\dim$
be any dimension function which assigns to a finitely generated left
$k*G$-module a nonnegative real number and satisfies
\begin{enumerate}
\item
$\dim(k*G) = 1$.
\item
For every short exact sequence $0 \to M_0 \to M_1 \to M_2 \to 0$
of finitely generated left $k*G$-modules, we get
\[
\dim(M_1) = \dim(M_0) + \dim(M_2).
\]
\item
If the finitely generated left $k*G$-module $M$ satisfies $\dim(M) =
0$, then every quotient module $Q$ of $M$ satisfies $\dim(Q) = 0$.
\end{enumerate}

Then for every finitely presented left $k*G$-module $M$, we get
$\dim(M) = \dim_{k*G}^{\Ore}(M)$.
\end{theorem}

\begin{proof}
Let $S$ denote the non-zero-divisors of $k*G$.
We have to show that for all $r,s \in \mathbb{N}$ and every $r\times
s$ matrix $A$ with entries in $k*G$
\begin{equation}
\dim_{k*G}^{\Ore}\bigl(\coker\big(r_A \colon S^{-1}k*G^r \to
S^{-1}k*G^s\bigr) \bigr)=
\dim\bigl(\coker\big(r_A \colon k*G^r \to k*G^s\bigr) \bigr),
\label{dim_Skg_is_dim}
\end{equation}
where $r_A$ denotes the module homomorphism given by right
multiplication with $A$.
First note that we may assume that $r=s$.  Indeed if $r<s$, replace
$A$ with the $s \times s$ matrix which is $A$ for the first $r$ rows,
and has 0's on the bottom $s-r$ rows.  On the other hand if $r>s$,
replace $A$ with the $r \times r$ matrix $B$ with entries $(b_{ij})$
which is $A$ for the first $s$ columns, and has $b_{ij} =
\delta_{ij}$ if $i>s$, where $\delta_{ij}$ is the Kronecker delta.

We will often use the obvious long exact sequence
associated to homomorphisms $f \colon M_0 \to M_1$ and $g \colon
M_1 \to M_2$
\begin{equation}
0 \to \ker(f) \to \ker(g\circ f) \to \ker(g) 
\to \coker(f) \to \coker(g\circ f) \to \coker(g) \to 0.
\label{ker-coker-sequence}
\end{equation}
We now assume that $A$ is an $r \times r$ matrix.  Note that
equation \eqref{dim_Skg_is_dim} is true if $A$ is invertible over
$S^{-1}k*G$; this is because then $\ker r_A = 0$ (whether $A$ is
considered as a matrix over $k*G$ or $S^{-1}k*G$).

Next observe that if $U \in \Mat_r(k*G)$ which is invertible over
$\Mat_r(S^{-1}k*G)$, then equation \eqref{dim_Skg_is_dim} holds for
$A$ if and only if it holds for $AU$, and also if and only if it hold
for $UA$.  This follows from
\eqref{ker-coker-sequence}, $\ker U = 0$, $\dim (\coker U) =
\dim_{k*G}^{\Ore} (\coker U) = 0$, and in the second case we use the
third property of $\dim$.

We may write $S^{-1}k*G = \Mat_d(D)$ for some $d\in \mathbb{N}$ and
some skew field $D$.
By applying the Morita equivalence from $\Mat_d(D)$ to $D$ and back and 
doing Gaussian elimination over $D$ we see that there are 
invertible
matrices $U,V \in \Mat_{rd}(S^{-1}k*G)$ such that
$U\diag(A,\dots,A)V = J$, where there are $d$ $A$'s and $J$ is a matrix of
the form $\diag(1,\dots,1,0,\dots,0)$.
Now choose $u,v \in S$ such that $uU,vV \in
\Mat_{rd}(k*G)$.  Then $(uU) \diag(A,\dots,A) (Vv) = uJv$,
and the result follows.
\end{proof}

\begin{theorem}[Comparing Elek's dimension and the Ore dimension]
  \label{the:Comparing_Eleks_dimension_and_the_Ore_dimension}
  Let $G$ be a finitely generated group and let $k$ be a
skew field. Suppose that $kG$ is a prime left Goldie ring.
  Then for any finitely presented left $kG$-module $M$
$$\dim^{\Elek}_{kG} (M) = \dim_{kG}^{\Ore}(M).$$
\end{theorem}
\begin{proof}
This follows from Theorem~\ref{the:Comparing_dimensions} and
Theorem~\ref{the:Eleks-dimension}.
\end{proof}


\typeout{--------------------------------- Section 4  ----------------------------------------}

\section{Proof of the main theorem}
\label{sec:Proof_of_the_main_theorem}

\begin{proof}[Proof of Theorem~\ref{the:dim_approximation_over_fields}]~%
\ref{the:dim_approximation_over_fields:finitely_presented_modules}
In the first step we reduce the claim to the case, where $G$ is finitely generated.
Consider a finitely presented left $kG$-module $M$. Choose a matrix $A \in M_{r,s}(kG)$
such that $M$ is $kG$-isomorphic to the cokernel of $r_A \colon kG^r \to kG^s$.
Since $A$ is a finite matrix and each element in $kG$ has finite support, we can find a
finitely generated subgroup $H \subseteq G$ such that $A \in
M_{r,s}(kH)$. Both $kG$ and $kH$ are 
prime left Goldie by Lemma~\ref{lem:delta_implies_prime} and 
Theorem~\ref{thm:elementary_amenable_goldie}.
Consider the finitely presented $kH$-module $N := \coker\bigl(r_A \colon kH^r \to kH^s\bigr)$.
Then $M = kG \otimes_{kH} N$. We can also consider the Ore localization
$T^{-1}kH$ for $T$ the set of non-zero-divisors of $kH$.
Put $H_n = H \cap G_n$. We obtain a  
residual chain $(H_n)_{n\ge 0}$ of $H$ and have: 
\begin{align*}
\dim_{kG}^{\Ore}(M)
&= 
\dim_{S^{-1}kG} \bigl(S^{-1}kG \otimes_{kG} M\bigr)\\
&= 
\dim_{S^{-1}kG} \bigl(S^{-1}kG \otimes_{kG} kG \otimes_{kH} N\bigr)
\\
&= 
\dim_{S^{-1}kG}\bigl(S^{-1}kG \otimes_{T^{-1}kH} T^{-1}kH \otimes_{kH} N\bigr)
\\
&= 
\dim_{T^{-1}kH}\bigl(T^{-1}kH \otimes_{kH}  N)
\\
&= 
\dim_{kH}^{\Ore}(N).
\end{align*}
We compute
\begin{align*}
\frac{\dim_k\bigl(k\otimes_{kG_n} M\bigr)}{[G:G_n]} 
&=  
\frac{\dim_k\bigl(k[G/G_n]\otimes_{kG} M\bigr)}{[G:G_n]}\\
&=  
\frac{\dim_k\bigl(k[G/G_n]\otimes_{kG} kG \otimes_{kH} N\bigr)}{[G:G_n]}
\\
&=  
\frac{\dim_k\bigl(k[G/G_n]\otimes_{k[H/H_n]} k[H/H_n] \otimes_{kH} N\bigr)}{[G:G_n]}
\\
&=  
\frac{[G/G_n:H/H_n] \cdot \dim_k\bigl(k[H/H_n] \otimes_{kH} N\bigr)}{[G:G_n]}
\\
&=  
\frac{[G/G_n:H/H_n] \cdot \dim_k\bigl(k \otimes_{kH_n} N\bigr)}{[G/G_n:H/H_n] \cdot [H:H_n]}
\\
&=  
\frac{\dim_k\bigl(k \otimes_{kH_n} N\bigr)}{[H:H_n]}.
\end{align*}
Therefore the claim holds for $M$ over $kG$ if it holds for $N$ over $kH$.
Hence we can assume without loss of generality that $G$ is finitely generated.

Now apply Theorem~\ref{the:Approximation_for_finitely_presented_kG-modules_for_Eleks_dimension_function}
and Theorem~\ref{the:Comparing_Eleks_dimension_and_the_Ore_dimension}.
\\[1mm]~\ref{the:dim_approximation_over_fields:chain_complexes}
We obtain from additivity, the exactness of the functor $S^{-1}kG \otimes_{kG}\_$
and the right exactness of the functor $k \otimes_{kG}\_$ that 
\begin{align*}
\dim^{\Ore}_{kG}\bigl(H_i(C_*)\bigr)
&= 
\dim^{\Ore}_{kG} \bigl(\coker(c_{i+1})\bigr)
+\dim^{\Ore}_{kG}\bigl(\coker(c_i)\bigr)
- \dim^{\Ore}_{kG}\bigl(C_{i-1}\bigr),\\
\dim_{k}\bigl(H_i(k \otimes_{kG_n} C_*)\bigr)
&= 
\dim_{k}\bigl(k \otimes_{kG_n}   \coker(c_{i+1})\bigr)
+\dim_{k}\bigl(k \otimes_{kG_n}  \coker(c_i)\bigr)
\\  & \hspace{10mm} - \dim_{k}\bigl(k \otimes_{kG_n} C_{i-1}\bigr).
\end{align*}
Hence the claim follows from assertion~\ref{the:dim_approximation_over_fields:finitely_presented_modules}
applied to the finitely presented $kG$-modules $\coker(c_{i+1})$, $\coker(c_i)$ and $C_{i-1}$.
\\[1mm]~\ref{the:dim_approximation_over_fields:CW-complexes}
This follows from assertion~\ref{the:dim_approximation_over_fields:chain_complexes} applied
to the cellular chain complex of~$X$.
\end{proof}


\typeout{--------------------------------- Section 5  ----------------------------------------}

\section{Extension to the  virtually torsionfree case}
\label{sec:Extension_to_the_virtually_torsionfree_case}

Next we explain how Theorem~\ref{the:dim_approximation_over_fields} can be
extended to the virtually torsionfree case.

For the remainder of this section let $k$ be a skew field,
let $G$ be an amenable group which possesses a
subgroup $H$ of finite index with $\Delta^+(H) = 1$,
and let $k*G$ be a crossed
product such that $k*H$ is a left Goldie ring.
We define the \emph{virtual Ore dimension} of a
$k*G$-module~$M$ by
\begin{equation}\label{eq:virtual_Ore_dimension}
  \vdim^{\Ore}_{k*G}(M) =
  \frac{\dim^{\Ore}_{k*H}\bigl(\res_{k*G}^{k*H} M\bigr)}{[G:H]},
\end{equation}
where $\res_{k*G}^{k*H} M$ is the $k*H$-module obtained from the $k*G$-module $M$ by
restricting the $G$-action to $H$. 

We have to show that this is independent of the choice of
$H$.  Since every subgroup of finite index contains a normal subgroup
of finite index, it is enough to show that if $K$ is a normal subgroup
of finite index in $H$ and $K \le H \le G$ with $H$ torsion free,
then for every $k*H$-module $N$,
\begin{equation}
 \frac{\dim^{\Ore}_{k*K}\bigl(\res_{k*H}^{k*K} N\bigr)}{[H:K]}
=
\dim^{\Ore}_{k*H}(N).
\label{Ore-dimension_and_restriction}
\end{equation}
Let $T$ denote the set of non-zero-divisors of $k*H$ and write
$S = (k*K) \cap T$.  Note that $\Delta^+(K) = 1$, so $k*K$ is still
a prime left Goldie ring and hence the ring $S^{-1}k*K$ exists.
Then $S^{-1}k*H \cong (S^{-1}k*K)*[H/K]$ and there is a natural ring
monomorphism $\theta\colon S^{-1}k*H \hookrightarrow T^{-1}k*H$.
Since $S^{-1}k*K$ is a matrix ring over a 
skew field by Theorem~\ref{thm:prime_Goldie_matrix},
we see that $(S^{-1}k*K)[H/K]$ is an Artinian ring,
because $H/K$ is finite.  But every element of $T$ is a
non-zero-divisor in $(S^{-1}k*K)[H/K]$, and since every
non-zero-divisor in an Artinian ring is 
invertible (compare~\cite{Rowen(2008)}*{Exercise~22 of Chapter~15 on p.~16}), 
we see that every
element of $T$ is invertible in $S^{-1}k*H$ and we conclude that
$\theta$ is onto and hence is an isomorphism.
We deduce that $\dim_{S^{-1}k*K}(T^{-1}k*H) = [H:K]$ and that the
natural map $S^{-1}N \to T^{-1}N$ induced by $s^{-1}n \mapsto
s^{-1}n$ is an isomorphism.  This proves
\eqref{Ore-dimension_and_restriction}.

\begin{theorem}[Extension to the virtually torsionfree case]
\label{the:Extension_to_the_virtually_torsionfree_case}
Let $G$ be an amenable group which possesses a
subgroup $E \subseteq G$ of finite index such that $kE$ is left
Goldie and $\Delta^+(E)=1$, and let $k$ be a skew field.
Then assertions~\ref{the:dim_approximation_over_fields:finitely_presented_modules},%
~\ref{the:dim_approximation_over_fields:chain_complexes}
and~\ref{the:dim_approximation_over_fields:CW-complexes} of
Theorem~\ref{the:dim_approximation_over_fields} remain true,
provided we replace $\dim_{kG}^{\Ore}$ by $\vdim_{kG}^{\Ore}$ everywhere.
\end{theorem}
\begin{proof}
It suffices to prove the claim for assertion~\ref{the:dim_approximation_over_fields:finitely_presented_modules}
since the proof in Theorem~\ref{the:dim_approximation_over_fields}
that it implies the other two assertions applies also in this more
general situation. Let $(G_n)_{n\ge 0}$ be a residual chain of $G$. 
To prove the result in general, we may assume that $G$ is finitely
generated.
Since $kE$ is left Goldie and $[G:E]<\infty$, the group ring 
$kG$ is also left Goldie. Further, every $kG_n$ is left Goldie. 
Since $\Delta^+(G)$ is finite (its order is bounded by $[G:E]$), 
there exists $N \in\mathbb{N}$ such that $G_N \cap \Delta^+(G) = 1$, and then
$\Delta^+(G_N) = 1$ and $G_i \subseteq G_N$ for all $i \ge N$.  Set
$H = G_N$, so that $kH$ is prime by Lemma~\ref{lem:delta_implies_prime}. 
Then for a finitely presented $kH$-module $L$,
\begin{equation*}
\dim_{kH}^{\Ore}(L) = 
\lim_{n \to \infty} \frac{\dim_k\bigl(k \otimes_{k[G_n\cap H]}
L\bigr)}{[H:H \cap G_n]}
\end{equation*}
by 
Theorems~\ref{the:Approximation_for_finitely_presented_kG-modules_for_Eleks_dimension_function} and~\ref{the:Comparing_dimensions}. 
We have $[G:G_n] = [G:H] \cdot [H : H \cap G_n]$ for $n \ge N$.
This implies for every finitely presented $kG$-module $M$
\begin{align*}
\vdim_{kG}^{\Ore}(M)
  = 
\frac{\dim_{kH}^{\Ore}\bigl(\res_{kG}^{kH} M\bigr)}{[G:H]}
& =  
\lim_{n \to \infty} \frac{\dim_k\bigl(k \otimes_{k[H\cap G_n]}
\res_{kG}^{kH}  M\bigr)}{[G:H] \cdot [H:H \cap G_n]}
\\
& =  
\lim_{n \to \infty} \frac{\dim_k\bigl(k \otimes_{kG_n}
M\bigr)}{[G:G_n]}.\qedhere
\end{align*}
\end{proof}

\begin{remark}
 Because of Theorem~\ref{thm:elementary_amenable_goldie},
Theorem~\ref{the:dim_approximation_over_fields} is true
in the case $k$ is a skew field and
$G$ is an elementary amenable group in which the orders of the finite 
subgroups are bounded (clearly $\Delta^+(G_n) = 1$ for sufficiently
large $n$).
In particular Theorem~\ref{the:dim_approximation_over_fields} 
is true for any virtually torsionfree elementary amenable group.

\end{remark}


\typeout{--------------------------------- Section 6  ----------------------------------------}

\section{Examples}
\label{sec:Examples}

\begin{remark}\label{rem:connection_to_L2_for_group_C}
Let $(G_n)_{n\ge 0}$ be a residual chain of a group $G$. 
Let $X$ be a finite free $G$-$CW$-complex.
Let $k$ be a field of characteristic $\charac(k)$. For a prime $p$ denote 
by $\IF_p$ the field of $p$ elements. 
Then we conclude from the universal coefficient theorem
\begin{align*}
\dim_{k}\bigl(H_i(G_n\backslash X;k)\bigr) & =  \dim_{\IQ}\bigl(H_i(G_n\backslash X;\IQ)\bigr)  &  \charac(k) = 0;
\\
\dim_{k}\bigl(H_i(G_n\backslash X;k)\bigr) & =  \dim_{\IF_p}\bigl(H_i(G_n\backslash X;\IF_p)\bigr)  &  p = \charac(k) \not= 0;
\\
\dim_{\IF_p}\bigl(H_i(G_n\backslash X;\IF_p)\bigr) & \ge  \dim_{\IQ}\bigl(H_i(G_n\backslash X;\IQ)\bigr).& 
\end{align*}
In particular we conclude from Remark~\ref{rem:fields_of_characteristic_zero} that 
\begin{equation*}
\liminf_{n \to \infty} \frac{\dim_{k}\bigl(H_i(G_n\backslash X;k)\bigr)}{[G:G_n]} 
\ge 
\lim_{n \to \infty} \frac{\dim_{k}\bigl(H_i(G_n\backslash X;\IQ)\bigr)}{[G:G_n]} 
= b_i^{(2)}(X;\caln(G)), 
\end{equation*}
where the latter term denotes the $i$-th $L^2$-Betti number of $X$. 
In particular we get from
Theorem~\ref{the:dim_approximation_over_fields}
for a torsionfree amenable group $G$ with 
no zero-divisors in $kG$ that 
\begin{align*}
\dim^{\Ore}_{kG}\bigl(H_i(X;k)\bigr) 
 = 
\lim_{n \to \infty} \frac{\dim_{k}(H_i(G_n\backslash X;k))}{[G:G_n]} 
& \ge 
\lim_{n \to \infty} \frac{\dim_{k}(H_i(G_n\backslash X;\IQ))}{[G:G_n]} 
\\
& = 
b_i^{(2)}(X;\caln(G))
\\
& = 
\dim^{\Ore}_{\IC G}\bigl(H_i(X;\IC)\bigr).
\end{align*}
This inequality is in general not an equality as the next example shows.
\end{remark}

\begin{example}
\label{exa:S1_wedge_Sd_cupDd_plus_1}
Fix an integer $d \ge 2$ and a prime number $p$.  Let $f_p \colon S^d \to S^d$
be a map of degree $p$ and denote by $i \colon S^d \to S^1 \vee S^d$ the
obvious inclusion.  Let $X$ be the finite $CW$-complex obtained from $S^1 \vee
S^d$ by attaching a $(d+1)$-cell with attaching map $i \circ f^d \colon S^d
\to S^1 \vee S^d$. Then $\pi_1(X) = \IZ$.  Let
$\widetilde{X}$ be the universal covering of $X$ which is a finite free
 $\IZ$-$CW$-complex.  Denote by $X_n$ the covering of $X$
associated to $n \cdot \IZ \subseteq \IZ$. The cellular $\IZ C$-chain complex of $\widetilde{X}$ is
concentrated in dimension $(d+1)$, $d$ and $1$ and $0$, the $(d+1)$-th
differential is multiplication with $p$ and the first differential is
multiplication with $(z-1)$ for a generator $z \in \IZ$
$$0\to \cdots \to \IZ[\IZ] \xrightarrow{p} \IZ[\IZ]  \to \cdots \to \IZ[\IZ]  \xrightarrow{z-1} \IZ[\IZ].$$
If the characteristic of $k$ is different from $p$, one easily checks that
$H_i(C_*) = 0$ and 
\[\dim^{\Ore}_{k \IZ}\bigl(H_i(\widetilde{X};k)\bigr) = 0~\text{ for  
$i \in \{d,d+1\}$}.\] 
If $p$ is the characteristic of $k$, then $H_i(C_*) = k\IZ$ and
\[\dim^{\Ore}_{k\IZ}\bigl(H_i(\widetilde{X};k)\bigr) = 1~\text{ for $i \in \{d,d+1\}$.}\]
Hence $\dim^{\Ore}_{kG}\bigl(H_i(\widetilde{X};k)\bigr)$ does depend on $k$ in
general.
\end{example}

\begin{example}\label{exa:s1_actions}
 Let $G$ be a torsionfree amenable group such that $kG$ has no
 zero-divisors. Let $S^1\to X\to B$ be a fibration  of connected $CW$-complexes such that 
$X$ has fundamental group $\pi_1(X)\cong G$ and  $\pi_1(S^1) \to \pi_1(X)$ is injective. Then 
\begin{equation}\label{eq:fiber_bundles}
\dim^{\Ore}_{kG}\bigl(H_i(\widetilde X;k)\bigr)=0
\end{equation}
for every $i\ge 0$. 

Let $S=kG-\{0\}$ and $S_0=k\IZ-\{0\}$. 
 By looking at the cellular chain complex one directly sees that 
	\[
		H_i\bigl(\widetilde S^1, S_0^{-1}k\IZ\bigr)=0~\forall
                i\ge 0, 
	\]
	thus $H_i\bigl(\widetilde S^1,S^{-1}kG\bigr)=S^{-1}kG\otimes_{S_0^{-1}k\IZ}H_i\bigl(\widetilde S^1,S_0^{-1}k\IZ\bigr) =0$ for every $i\ge 0$. 
	The assertion is implied by the Hochschild-Serre spectral sequence that 
	converges to $H_{p+q}(\widetilde X,S^{-1}kG)$ and has the $E^2$-term: 
	\[
		E^2_{pq}=H_p\bigl(\widetilde B, H_q(\widetilde S^1, S^{-1}kG)\bigr).
	\]
\end{example}

\begin{example}[Sublinear growth of Betti numbers] 
Let $G$ be an infinite amenable group which possesses a
subgroup $H$ of finite index such that $kH$ is left Goldie 
and $\Delta^+(H) = 1$, e.g., 
$G$ is a virtually torsionfree elementary amenable group.
Let $k$ be a field. Let $(G_n)_{n\ge 0}$ be a residual 
chain of $G$. 
Denote by $b_i(G/G_n;K)$ the $i$-th Betti number of
the group $G/G_n$ with coefficients in $k$. 
Then we get for every $i \ge 0$
\begin{equation*}
\lim_{n \to \infty} \frac{b_i(G/G_n;k)}{[G:G_n]} = 0.
\end{equation*}
For $i=0$ this is obvious. For $i\ge 1$ this  follows 
from Theorem~\ref{the:Extension_to_the_virtually_torsionfree_case} and 
$H_i(EH;k)=H_i(H;k)= 0$. 
\end{example}


\typeout{-------------------------------- References  ---------------------------------------}



\begin{bibdiv}
\begin{biblist}

\bib{Abert-Jaikin-Zapirain-Nikolov(2007)}{unpublished}{
      author={Abert, M.},
      author={Jaikin-Zapirain, A.},
      author={Nikolov, N.},
       title={The rank gradient from a combinatorial viewpoint},
        date={2007},
        note={arXiv:math/0701925v2},
}

\bib{Deninger-Schmidt(2007)}{article}{
      author={Deninger, Christopher},
      author={Schmidt, Klaus},
       title={Expansive algebraic actions of discrete residually finite
  amenable groups and their entropy},
        date={2007},
        ISSN={0143-3857},
     journal={Ergodic Theory Dynam. Systems},
      volume={27},
      number={3},
       pages={769\ndash 786},
         url={http://dx.doi.org/10.1017/S0143385706000939},
      review={\MR{MR2322178 (2008d:37009)}},
}

\bib{Dodziuk-Linnell-Mathai-Schick_Yates(2003)}{article}{
      author={Dodziuk, J{\'o}zef},
      author={Linnell, Peter},
      author={Mathai, Varghese},
      author={Schick, Thomas},
      author={Yates, Stuart},
       title={Approximating {$L^2$}-invariants and the {A}tiyah conjecture},
        date={2003},
        ISSN={0010-3640},
     journal={Comm. Pure Appl. Math.},
      volume={56},
      number={7},
       pages={839\ndash 873},
         url={http://dx.doi.org/10.1002/cpa.10076},
        note={Dedicated to the memory of J{\"u}rgen K. Moser},
      review={\MR{MR1990479 (2004g:58040)}},
}

\bib{Elek(2002)}{article}{
      author={Elek, G{\'a}bor},
       title={Amenable groups, topological entropy and {B}etti numbers},
        date={2002},
        ISSN={0021-2172},
     journal={Israel J. Math.},
      volume={132},
       pages={315\ndash 335},
      review={\MR{MR1952628 (2003k:37026)}},
}

\bib{Elek(2003c)}{article}{
      author={Elek, G{\'a}bor},
       title={The rank of finitely generated modules over group algebras},
        date={2003},
        ISSN={0002-9939},
     journal={Proc. Amer. Math. Soc.},
      volume={131},
      number={11},
       pages={3477\ndash 3485 (electronic)},
      review={\MR{MR1991759 (2004i:43003)}},
}

\bib{Elek(2006strong)}{article}{
      author={Elek, G{\'a}bor},
       title={The strong approximation conjecture holds for amenable groups},
        date={2006},
        ISSN={0022-1236},
     journal={J. Funct. Anal.},
      volume={239},
      number={1},
       pages={345\ndash 355},
         url={http://dx.doi.org/10.1016/j.jfa.2005.12.016},
      review={\MR{MR2258227 (2007m:43001)}},
}

\bib{Kropholler-Linnell-Moody(1988)}{article}{
      author={Kropholler, P.~H.},
      author={Linnell, P.~A.},
      author={Moody, J.~A.},
       title={Applications of a new ${K}$-theoretic theorem to soluble group
  rings},
        date={1988},
        ISSN={0002-9939},
     journal={Proc. Amer. Math. Soc.},
      volume={104},
      number={3},
       pages={675\ndash 684},
      review={\MR{89j:16016}},
}

\bib{Lang(2002_Algebra)}{book}{
      author={Lang, Serge},
       title={Algebra},
     edition={third},
      series={Graduate Texts in Mathematics},
   publisher={Springer-Verlag},
     address={New York},
        date={2002},
      volume={211},
        ISBN={0-387-95385-X},
      review={\MR{MR1878556 (2003e:00003)}},
}

\bib{Lindenstrauss-Weiss(2000)}{article}{
      author={Lindenstrauss, Elon},
      author={Weiss, Benjamin},
       title={Mean topological dimension},
        date={2000},
     journal={Israel J. Math.},
      volume={115},
       pages={1\ndash 24},
      review={\MR{MR1749670 (2000m:37018)}},
}

\bib{Linnell(2006)}{incollection}{
      author={Linnell, Peter~A.},
       title={Noncommutative localization in group rings},
        date={2006},
   booktitle={Non-commutative localization in algebra and topology},
      series={London Math. Soc. Lecture Note Ser.},
      volume={330},
   publisher={Cambridge Univ. Press},
     address={Cambridge},
       pages={40\ndash 59},
      review={\MR{MR2222481 (2007f:16064)}},
}

\bib{Lueck(1994c)}{article}{
      author={L{\"u}ck, Wolfgang},
       title={Approximating ${L}\sp 2$-invariants by their finite-dimensional
  analogues},
        date={1994},
        ISSN={1016-443X},
     journal={Geom. Funct. Anal.},
      volume={4},
      number={4},
       pages={455\ndash 481},
      review={\MR{95g:58234}},
}

\bib{Lueck(2002)}{book}{
      author={L{\"u}ck, Wolfgang},
       title={{$L\sp 2$}-{I}nvariants: {T}heory and {A}pplications to
  {G}eometry and \mbox{{$K$}-{T}heory}},
      series={Ergebnisse der Mathematik und ihrer Grenzgebiete. 3.~Folge. A
  Series of Modern Surveys in Mathematics [Results in Mathematics and Related
  Areas. 3rd Series. A Series of Modern Surveys in Mathematics]},
   publisher={Springer-Verlag},
     address={Berlin},
        date={2002},
      volume={44},
        ISBN={3-540-43566-2},
      review={\MR{1 926 649}},
}

\bib{Pape(2008)}{article}{
      author={Pape, Daniel},
       title={A short proof of the approximation conjecture for amenable
  groups},
        date={2008},
        ISSN={0022-1236},
     journal={J. Funct. Anal.},
      volume={255},
      number={5},
       pages={1102\ndash 1106},
         url={http://dx.doi.org/10.1016/j.jfa.2008.05.018},
      review={\MR{MR2455493 (2010a:46171)}},
}

\bib{Passman(1977)}{book}{
      author={Passman, Donald~S.},
       title={The algebraic structure of group rings},
   publisher={Wiley-Interscience [John Wiley \&\ Sons]},
     address={New York},
        date={1977},
        ISBN={0-471-02272-1},
        note={Pure and Applied Mathematics},
      review={\MR{81d:16001}},
}

\bib{Passman(1986)}{book}{
      author={Passman, Donald~S.},
       title={Group rings, crossed products and {G}alois theory},
      series={CBMS Regional Conference Series in Mathematics},
   publisher={Published for the Conference Board of the Mathematical Sciences,
  Washington, DC},
        date={1986},
      volume={64},
        ISBN={0-8218-0714-5},
      review={\MR{MR840467 (87e:16033)}},
}

\bib{Rowen(2008)}{book}{
      author={Rowen, Louis~Halle},
       title={Graduate algebra: noncommutative view},
      series={Graduate Studies in Mathematics},
   publisher={American Mathematical Society},
     address={Providence, RI},
        date={2008},
      volume={91},
        ISBN={978-0-8218-0570-1},
      review={\MR{MR2462400 (2009k:16001)}},
}

\bib{Stenstroem(1975)}{book}{
      author={Stenstr{\"o}m, Bo},
       title={Rings of quotients},
   publisher={Springer-Verlag},
     address={New York},
        date={1975},
        note={Die Grundlehren der Mathematischen Wissenschaften, Band 217, An
  introduction to methods of ring theory},
      review={\MR{52 \#10782}},
}

\bib{Weiss(2001)}{incollection}{
      author={Weiss, Benjamin},
       title={Monotileable amenable groups},
        date={2001},
   booktitle={Topology, ergodic theory, real algebraic geometry},
      series={Amer. Math. Soc. Transl. Ser. 2},
      volume={202},
   publisher={Amer. Math. Soc.},
     address={Providence, RI},
       pages={257\ndash 262},
      review={\MR{MR1819193 (2001m:22014)}},
}

\end{biblist}
\end{bibdiv}

\end{document}